\documentclass[12pt]{article}
\usepackage{amssymb}
\usepackage{amsthm}
\usepackage{graphicx}
\usepackage{psfrag}
\usepackage{graphics}
\usepackage{amsthm}

\input pictex.tex

\newcommand{\ce}{\mathrm{C}}
\newcommand{\esp}{\hspace{0.06cm}}

\theoremstyle{definition}

\begin{document}

\date{}
\author{Andr\'es Navas}

\title{Wandering domains for diffeomorphisms of the $k$-torus:\\ 
a remark on a theorem by Norton and Sullivan}

\maketitle

\vspace{-0.2cm}

\noindent{\bf Abstract:} We show that there is no $C^{k+1}$ diffeomorphism 
of the $k$-torus which is semiconjugate to a minimal translation and has a 
wandering domain all of whose iterates are Euclidean balls.

\vspace{0.2cm}


\vspace{0.2cm}


\vspace{0.6cm}

\noindent{\Large{\bf Introduction}}

\vspace{0.5cm}

Answering (in the negative) a question raised by Poincar\'e \cite{Po}, Denjoy proved one of his most 
famous theorems \cite{De}, namely the action of every $C^2$ (orientation preserving) 
diffeomorphism of the circle with irrational rotation number is minimal. 
This result can be considered as one of the the starting points of the theory 
of differentiable dynamics, and many generalizations have 
been proposed in the context of one-dimensional dynamics (see, for example,  
\cite{demelo,DKN,ghys,Sa,chuarts,yoccoz}). The search for higher dimensional 
analogues of the Denjoy Theorem is a natural problem that has attracted some 
interest in recent years. Although there is some partial evidence in the positive 
direction, no definitive result is known. 

Let us be more precise. As is well known, the $C^2$ (or $C^{1+Lip}$) hypothesis 
is necessary in the statement of the 
Denjoy Theorem. Indeed, Denjoy himself gave $C^1$ 
counter-examples for his result. (It should be noticed that the first $C^1$ 
counter-examples were constructed by Bohl in \cite{bohl}.) These examples were 
improved by Herman in \cite{herman} up to the class $C^{2-\varepsilon}$ for every 
$\varepsilon > 0$ (se also \cite{tsuboi}), thus showing that the Denjoy Theorem is 
sharp in the H\"older scale. The blowing-up method of Denjoy and Herman is 
classical and direct. However, circle diffeomorphisms with similar properties 
can be obtained as holonomy maps along stable manifolds of Anosov 
diffeomorphisms of the 2-torus. Mac Swiggen extended this construction in 
\cite{MaS1,MaS2} to higher dimensions and showed that, for each $k \geq 2$,  
there exists a dense family $(\theta_1,\ldots,\theta_k)$ of 
$\mathbb{Q} / \mathbb{Z}$-independent numbers in 
$\mathbb{R} / \mathbb{Z}$ such that, for any $\varepsilon \!>\! 0$, there 
is a $\ce^{k + 1 - \varepsilon}$ diffeomorphism of the $k$-torus that has a 
wandering (topological) disk and is semiconjugate to the translation by 
\esp $(\theta_1,\ldots,\theta_k)$. Unfortunately, the family of rotation vectors 
which appear in Mac Swiggen's constructions is countable (they are all algebraic, 
and therefore Diophantine), and it is unclear whether similar examples do exist 
for any translation vector with the above properties.

Mac Swiggen's examples show that the natural differentiability where we should 
look for an analogue of the Denjoy Theorem on the $k$-torus is $C^{k+1}$. This 
is confirmed by the fact that, by a straightforward application of the KAM Theory, 
if $f$ is a small $C^{k+1+\varepsilon}$ perturbation of an irrational 
Diophantine translation of the $k$-torus, then $f$ is $C^{k}$ conjugate to it. 

In the one-dimensional case, the blowing-up procedure for the construction 
of counter-examples is necessarily ``conformal'': one replaces points along 
orbits by intervals. In higher dimensions, one should replace points by continua, 
and the case where these continua are not topological disks is interesting by itself: 
see, for instance, \cite{sylvain,rees}. (We should point out that, in this case, 
perhaps a Denjoy type theorem holds in regularity smaller than $C^{k+1}$.) In 
the case of wandering topological disks, several partial results are known in  
dimension 2 (see for example \cite{norton,norton-sullivan,nortonw}). 
For instance, in \cite{norton-sullivan}, Norton and Sullivan show   
that it is impossible for a $C^3$ diffeomorphism of the 
$2$-torus to be semiconjugate and non-conjugate to a minimal 
translation, provided that the preimages of points have some ``uniform'' 
conformal geometry along the orbits.

The aim of this Note is to show how a simple modification of some of the ideas of \cite{norton-sullivan} 
allows proving similar results for $C^{k+1}$ diffeomorphisms of the $k$-torus. It should be emphasized that 
these do not follow from Norton-Sullivan's arguments, as these strongly rely on the Morrey-Bojarsky-Ahlfors-Bers 
integration theorem, which is no longer available in higher dimensions. In particular, one of the key arguments of  
\cite{norton-sullivan} uses Sullivan's integrability theorem, which states that every uniformly quasi-conformal group 
of homeomorphisms of a surface is quasiconformally conjugate to a group of conformal maps \cite{sullivan}, and which is 
known to be false in dimension greater than two  \cite{tukia}.

For the sake of concreteness, we only prove the following theorem, which is somewhat the core of \cite{norton-sullivan}.

\vspace{0.65cm}

\noindent{\bf Theorem.} {\em Let $k \geq 2$, and let $f$ be a diffeomorphism of the 
$k$-torus that is semiconjugate to a minimal translation without being conjugate 
to it. If the preimage by the semiconjugacy of each point is either a point or an  
Euclidean ball, then $f$ cannot be of class $C^{k+1}$.}

\vspace{0.65cm}

Note that this result is still true in dimension one (where it follows from the 
classical Denjoy Theorem), but our arguments only work in higher dimensions 
(see, however, \cite[Exercise 3.1.4]{libro}, which is somewhat related to our arguments here). 
  
We should stress that we do not know whether the $C^{k+1}$ regularity hypothesis is actually needed for our Theorem 
(assuming that the wandering domains are Euclidean balls). Indeed, in Mac Swiggen's examples, the wandering domains 
have a very irregular geometry, and there is even no uniform bound for the diameter of their lifts. See also \cite{marko} for 
a recent interesting result concerning topological entropy of this kind of maps. 

Combining the methods of this Note with those of \cite{rusia}, one can show that if $f_i$, 
$i \!\in\! \{1,\ldots,d\}$, are respectively $C^{1+k_i}$ diffeomorphisms of the $k$-torus 
that are semiconjugate to minimal translations and whose translation vectors are 
independent over $\mathbb{Q} / \mathbb{Z}$, then one has $k_1 + \ldots + k_d \leq k$  
provided the $f_i$'s commute and the preimages of points by the (simultaneous) 
conjugacy to translations are either points or Euclidean balls. (The constants $k_i$'s 
are supposed to be positive but not necessarily integer numbers). Actually, the same statement 
holds without the hypothesis that the $f_i$ commute but asking for the commutativity of their 
permutation action along an orbit of balls arising from the blowing-up procedure.

\vspace{0.6cm}

\noindent{\bf Acknowledgments.} This Note circulated as a manuscript more than ten years ago. I'm indebted to 
J. Kiwi, V. Kleptsyn and M. Ponce for their comments at that time, to all my colleges who asked me to make this 
available (despite no progress has been made since then), and to J.~Bochi for his insight concerning Lemma 1.

The preparation of this text was funded by the CONICYT Project 1415 ``Geometr\'{\i}a en La Frontera'' and the 
Fondecyt grant 11060541. I also whish to thank the Institute of Pure Mathematics of Teheran 
(Iran) for the hospitality during this task, and M. Nassiri for his invitation.


\vspace{0.8cm}

\noindent{\Large {\bf Proof of the Theorem}}

\vspace{0.4cm}

Let $\mathrm{Conf} (k) \sim \mathrm{GL} (k,\mathbb{R}) / 
(\mathrm{SO} (k,\mathbb{R}) \times \mathbb{R}) $ denote the space of conformal 
structures on $\mathbb{R}^k$. This is a simply-connected space that carries a nonpositively curved metric which is invariant under 
the $\mathrm{GL} (k,\mathbb{R})$-action given by $A \cdot [B] := [A B]$. In particular, the distance function \esp $\mathrm{dist}_k$ \esp 
on it is smooth. 
For simplicity, we denote $\sigma_0 := [Id]$.  

\vspace{0.65cm}

\noindent{\bf Lemma 1.} {\em If $f$ is a $C^k$ diffeomorphism satisfying the hypothesis of the Theorem, then there exists 
a constant $M$ such that \esp $\mathrm{dist}_k ( [Df^n (x)], \sigma_0 ) \leq M$ \esp holds for all $x \in \mathbb{T}^k$ and 
all $n \geq 1$.}

\vspace{0.65cm}

\noindent{\em Proof.} 
Let $\varphi$ denote the semiconjugacy of $f$ to the corresponding translation. Then  
$$\Gamma := \mathbb{T}^k \setminus \bigcup \mathrm{interior} 
\big( \{ \varphi^{-1}(x) \!: x \in \mathbb{T}^k \} \big)$$
is a connected, nonwhere dense, minimal invariant set for 
$f$  (see \cite{nortonw}). Moreover, since all the wandering topological balls 
are Euclidean, $[Df(x)]$ is identically equal to $\sigma_0$ on $\Gamma$. Furthermore, 
as $\mathrm{dist}_k$ is smooth, the derivatives of the function \esp 
$x \mapsto \mathrm{dist}_k \big( [Df(x)], \sigma_0 \big)$ \esp 
vanish up to order $k$. 
By a successive application of the Mean Value Theorem, this implies that there exists a constant $C$ such 
that, if $x$ belongs to the interior of a ball $B_x$ that collapses to a single point under $\varphi$, then 
$$\mathrm{dist}_k ([Df(x)] , \sigma_0)  \leq C \esp \ell(x)^k,$$
where $\ell(x)$ is the half of the length of the shortest chord of $B_x$ through $x$. 
Since \esp $\ell(x) \!\leq\! \mathrm{radius} (B_x)$, this shows that
$$\mathrm{dist}_k ([Df(x)] , \sigma_0)  \leq M \esp \mathrm{vol} (B_x)$$
for a certain constant $M$. This yields, for each $n \geq 1$ and every $x \notin \Gamma$,
\begin{eqnarray*}
\mathrm{dist}_k ([Df^n(x)], \sigma_0) 
&=& \mathrm{dist}_k (Df^{n}(x) \cdot \sigma_0 , \sigma_0) \\
&\leq& 
\sum_{i=0}^{n-1} \mathrm{dist}_k \big( Df^{i+1}( f^{n-i-1}(x)) \cdot \sigma_0, Df^{i}(f^{n-i}(x)) \cdot \sigma_0 \big) \\
&=& \sum_{i=0}^{n-1} \mathrm{dist}_k ( Df^{i}(f^{n-i}(x))) \cdot [ Df (f^{n-i-1}(x)] , Df^i (f^{n-i}(x)) \cdot \sigma_0 ) \\
&=& \sum_{i=0}^{n-1} \mathrm{dist}_k ( [Df (f^{n-i-1}(x)) ], \sigma_0 ) \\
&\leq& \sum_{i=0}^{n-1} M \esp \mathrm{vol} (B_{f^{n-i-1}(x)}) \\
&\leq& M,
\end{eqnarray*}
where the last inequality holds because the balls $B_{f^{j}(x)}$ are two-by-two disjoint. 
Since $\Gamma$ is nonwhere dense in $\mathbb{T}^k$, the estimate 
above holds for every $x \!\in\! \mathbb{T}^k$. 

\vspace{0.65cm}

\noindent{\bf Remark.} The {\em dilatation} of an invertible linear map $A: \mathbb{R}^k \to \mathbb{R}^k$ is defined as
$$dil (A) := \frac{\max_{\|v\|=1} \| A(v) \|}{\min_{\|w\|=1} \|A(w)\|}.$$
This induces a function on $\mathrm{Conf} (n)$ that measures the degree of non-conformality of matrices. 
In the 2-dimensional case, this is a smooth function, as is shown by the well-known formula
$$dil (A) = \frac{1 + \|\mu\|}{1-\|\mu\|} = \exp (\mathrm{dist}_{hyp} (0,\mu)),$$
where $\mu$ denotes the Beltrami differential and $\mathrm{dist}_{hyp}$ stands for the hyperbolic distance on the 
Poincar\'e disk. However, in higher dimension, this function is no longer smooth (yet it is 
locally Lipschitz). This is why, unlike \cite{norton-sullivan}, we do not deal with the 
function $dil$, and we directly consider the function \esp $\mathrm{dist}_k$.

\vspace{0.65cm}

\noindent{\bf Lemma 2.} {\em For each $M > 0$ there exists \esp $\lambda > 1$ \esp and 
\esp $1 < \lambda'$ \esp with the following property: if $g$ is a diffeomorphism of 
$\mathbb{R}^k$ that commutes with the translations by vectors in $\mathbb{Z}^k$ 
and such that \esp $\mathrm{dist}_k ([Dg(x)],\sigma_0) \leq M$ \esp for all \esp $x \in \mathbb{R}^k$ \esp 
and \esp $g( B(x_0,\alpha) ) = B(y_0,\beta)$ \esp for some \esp $x_0,y_0$ \esp 
in $\mathbb{R}^k$ \esp and some positive numbers $\alpha,\beta$, then} 
$$g (B (x_0,\lambda \alpha)) \subset B(y_0,\lambda' \beta).$$

\noindent{\em Proof.} This follows directly from the equicontinuity of the family of restrictions to 
$B(x_0, 2 \alpha)$ of the maps $g$ satisfying the properties above; see \cite[\S 19]{vaisala}.

\vspace{0.65cm}

To complete the proof of the Theorem, we will use an argument that, in the one-dimensional context, 
goes back to Schwartz \cite{chuarts}. Fix a wandering ball \esp $B \!=\! B(x_0,\alpha_0)$, \esp and 
denote by $x_n$ (resp. $\alpha_n$) the center (resp. the radius) of $f^n(B)$. Clearly, 
$\alpha_n$ goes to zero as $n$ goes to infinite, and for each $\varepsilon \!>\! 0$ one 
has $f^{n_k}(B) \!\subset\! B(x_0,\alpha_0 + \varepsilon)$ for an increasing sequence 
$(n_k)$ of positive numbers. Fix such an $n = n_k$ so that 
$$\alpha_{n} < \frac{(\lambda - 1) \alpha_0}{2 \lambda'} \qquad \mbox{ and } \qquad 
\mathrm{dist} (x_0,x_{n}) < \alpha_0 + \left( \frac{\lambda - 1}{2} \right) \alpha_0.$$ 
We can apply Lemma 2 to any covering map $\tilde{f}$ of $f$, thus obtaining 
$$\tilde{f}^{n_k} \Big( B(\tilde{x}_0,\lambda \alpha_0) \Big) 
\subset B \big( \tilde{x}_{n}, \lambda' \alpha_{n} \big) 
\subset B \left( \tilde{x}_0, \alpha_0 + \Big( \frac{\lambda - 1}{2} \Big) 
            \alpha_0 + \lambda' \alpha_{n} \right) 
\subset B (\tilde{x}_0, \lambda \alpha_0).$$
Therefore, 
$f^{n} (\overline{ B(x_0,\lambda \alpha_0)} ) \subset \overline{B (x_0, \lambda \alpha_0)}$, 
and by Brouwer's fixed point theorem, $f^{n}$ has a fixed point in 
$\overline{B (x_0, \lambda \alpha_0)}.$ However, this is absurd, since 
$f$ is semiconjugate to a minimal (and therefore periodic-point  
free) torus translation. This contradiction completes the proof.

\vspace{-0.5cm}

\begin{figure}[htp]
\hspace{-2.2cm}
\includegraphics[scale=0.8]{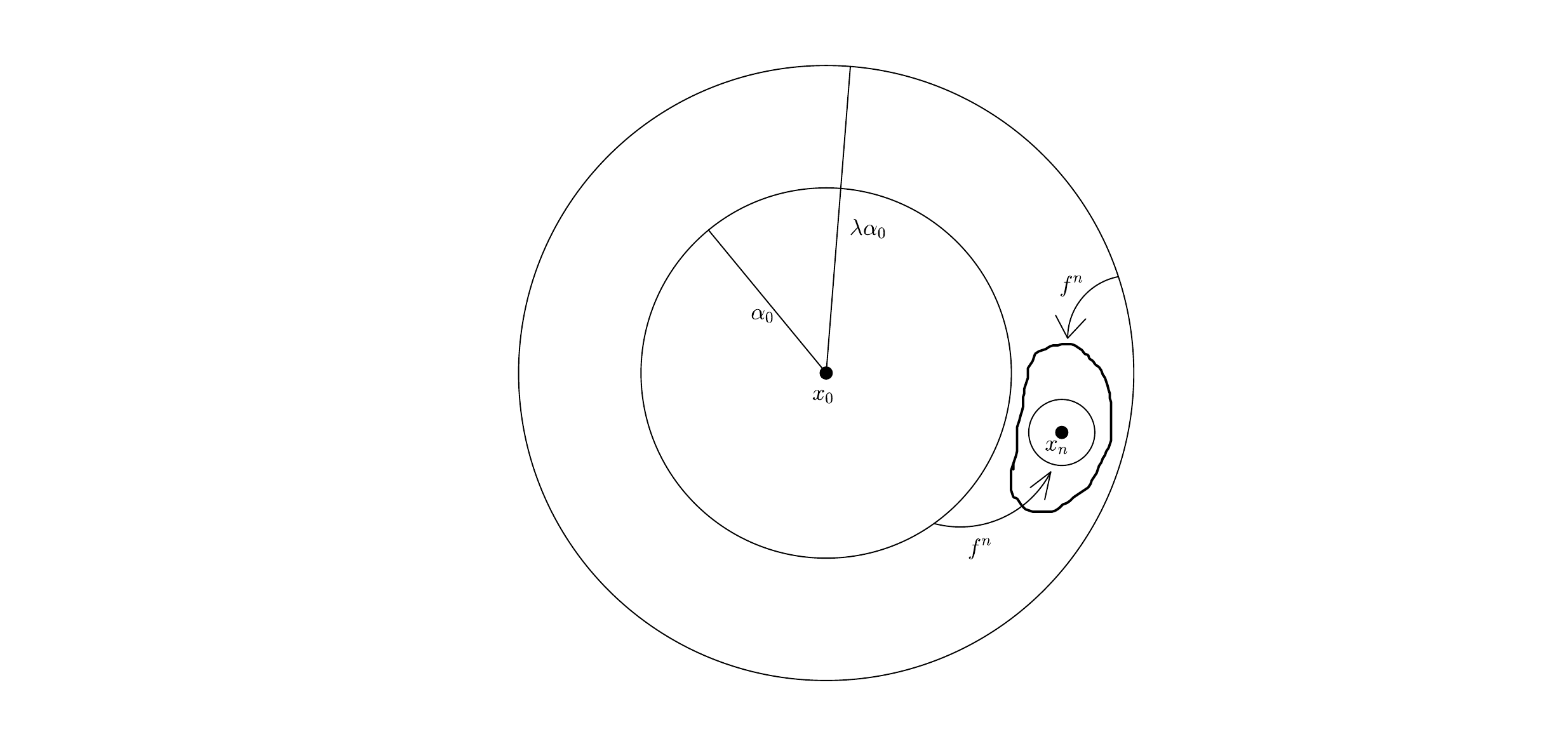}
\end{figure}


\vspace{-1cm}

\begin{footnotesize}

\vspace{0.1cm}

\noindent Andr\'es Navas\\

\noindent Dpto de Matem\'atica y C.C., Fac. de Ciencia, Univ. de Santiago de Chile\\ 

\noindent Alameda 3363, Estaci\'on Central, Santiago, Chile\\ 

\noindent Email address: andres.navas@usach.cl

\end{footnotesize}


\begin{thebibliography}{Dillo 83}

\bibitem{sylvain} {\sc B\'eguin, F., Crovisier, S., Le Roux, F.} Construction of 
curious minimal uniquely ergodic homeomorphisms on manifolds: the Denjoy-Rees 
technique. {\em Annales Scient. de l'\'Ecole Norm. Sup.} {\bf 40} (2007), 251-308.

\bibitem{bohl} {\sc Bohl, P.} Uber die hinsichtlich der unabh\"angigen variabeln
periodische. {\em Acta Math.} {\bf 40} (1916), 321-336.\\

\bibitem{demelo} {\sc De Melo, W. \& Van Strien, S.} A structure theorem in 
one dimensional dynamics. {\em Annals of Math.} {\bf 129} (1989), 519-546.\\

\bibitem{De} {\sc Denjoy, A.} Sur les courbes d\'efinies par des \'equations 
diff\'erentielles \`a la surface du tore. {\em J. Math. Pure Apl.} {\bf 11} 
(1932), 333-375.\\

\bibitem{DKN} {\sc Deroin, B., Kleptsyn, V. \& Navas, A.} Sur la 
dynamique unidimensionelle en r\'egularit\'e interm\'ediaire. 
{\em Acta Math.} {\bf 199} (2007), 199-262.\\

\bibitem{ghys} {\sc Ghys, \'E.} Transformations holomorphes au voisinage d'une 
courbe de Jordan. {\em Comptes Rendus Acad. Sci. Paris} {\bf 298} (1984), 385-388.\\

\bibitem{herman} {\sc Herman, M.} Sur la conjugaison diff\'erentiable des
diff\'eomorphismes du cercle \`a des rotations. {\em Publ. Math. de l'IH\'ES} 
{\bf 49} (1979), 5-234.\\

\bibitem{rusia} {\sc Kleptsyn, V. \& Navas, A.} A Denjoy type theorem for commuting 
circle diffeomorphisms with different H\"older differentiability classes. {\em Moscow 
Math. Journal} {\bf 8}, 477-492 (2008).\\

\bibitem{marko} {\sc Kwakkel, F. \& Markovic, V.} Topological entropy and 
diffeomorphisms of surfaces with wandering domains.  {\em Ann. Acad. Sci. Fenn. Math.} 
{\bf 35} (2010), no. {\bf 2}, 503-513.\\
  
\bibitem{MaS1} {\sc Mac Swiggen, P.} Diffeomorphisms of the torus with wandering
domains. {\em Proc. of the AMS} {\bf 117} (1993), 1175-1186.\\

\bibitem{MaS2} {\sc Mac Swiggen, P.} Diffeomorphisms of the $k$-torus with wandering
domains. {\em Erg. Theory and Dynam. Systems} {\bf 15} (1995), 1189-1205.\\ 

\bibitem{libro} {\sc Navas, A.} {\em Groups of circle diffeomorphisms.} 
Chicago Lectures in Mathematics (2011).\\ 

\bibitem{norton} {\sc Norton, A.} An area approach to wandering domains. 
{\em Erg. Theory and Dynam. Systems} {\bf 11} (1991), 455-467.\\

\bibitem{norton-sullivan} {\sc Norton, A. \& Sullivan, D.} Wandering domains and invariant 
conformal structures for mappings of the 2-torus. {\em Ann. Acad. Sc. Fenn.} {\bf 21} 
(1996), 51-68.\\

\bibitem{nortonw} {\sc Norton, A. \& Welling, J.} Conformal irregularity for
diffeomorphisms of the 2-torus. {\em Rocky Mountain J. Math.} {\bf 24} (1994), 651-671\\

\bibitem{Po} {\sc Poincar\'e, H.} M\'emoire sur les courbes d\'efinies par une 
\'equation diff\'erentielle III. {\em J. Math. Pure Apl.} {\bf 1} (1885), 167-244.\\

\bibitem{rees} {\sc Rees, M.} A minimal positive entropy homeomorphism of the 2-torus. 
{\em J. London Math. Soc.} {\bf 23} (1981), 537-550.

\bibitem{Sa} {\sc Sacksteder, R.} Foliations and pseudogroups. {\em Amer. J. Math.} 
{\bf 87} (1965), 79-102.\\

\bibitem{chuarts} {\sc Schwartz, A.} A generalization of Poincar\'e-Bendixon theorem to 
closed two dimensional manifolds. {\em Amer. J. Math.} {\bf 85} (1963), 453-458.\\

\bibitem{sullivan} {\sc Sullivan, D.} Riemann surfaces and related topics, in Proc. Stony Brook Conf., Stony Brook, 
N.Y.  (1978). {\em Ann. of Math. Studies} {\bf 97}, Princeton Univ. Press, Princeton, N.J. (1981), 465-496.


\bibitem{tsuboi} {\sc Tsuboi, T.} Homological and dynamical study on certain groups of
Lipschitz homeomorphisms of the circle. {\em J. Math. Soc. Japan} {\bf 47} (1995), 1-30.\\

\bibitem{tukia} {\sc Tukia, P.} A quasiconformal group not isomorphic to a M\"obius group. 
{\em Ann. Acad. Sci. Fenn. Ser. A I Math.} {\bf 6} (1981), no. {\bf 1}, 149-160.\\

\bibitem{vaisala} {\sc V\"ais\"al\"a, J.} Lectures on $n$-dimensional quasiconformal 
mappings. {\em Lecture Notes in Math.} {\bf 229} (1971).

\bibitem{yoccoz} {\sc Yoccoz, J-C.} Il n'y a pas de contre-exemple de Denjoy
analytique. {\em Comptes Rendus Acad. Sci. Paris} {\bf 298} (1984), 141-144.\\

\end{thebibliography}
\end{document}